\documentclass[12pt]{article}

\makeatletter
\renewcommand{\@oddfoot}{\hfill \thepage}
\makeatother

\usepackage[T2A]{fontenc}
\usepackage[cp1251]{inputenc}
\usepackage[russian, english]{babel}
\usepackage{mathrsfs}
\usepackage{amssymb}
\usepackage{amsmath,amsthm,amsfonts}
\usepackage[dvipdf]{graphicx}
\usepackage{fancyhdr}

\begin{document}

\begin{center}
{\bf PROBABILITY DISTRIBUTION FUNCTION\\ FOR THE EUCLIDEAN DISTANCE \\ BETWEEN TWO TELEGRAPH PROCESSES}
\end{center}

\begin{center}
Alexander D. KOLESNIK\\
Institute of Mathematics and Computer Science\\
Academy Street 5, Kishinev 2028, Moldova\\
E-Mail: kolesnik@math.md
\end{center}

\vskip 0.2cm

\begin{abstract}
Consider two independent Goldstein-Kac telegraph processes $X_1(t)$ and $X_2(t)$ on the real line $\Bbb R$.
The processes $X_k(t), \; k=1,2,$ are performed by stochastic motions at finite constant velocities $c_1>0, \; c_2>0,$
that start at the initial time instant $t=0$ from the origin of the real line $\Bbb R$ and are controlled by two independent homogeneous Poisson processes 
of rates $\lambda_1>0, \; \lambda_2>0$, respectively. Closed-form expression for the probability distribution function of the Euclidean distance
$$\rho(t) = \vert X_1(t) - X_2(t) \vert , \qquad t>0,$$
between these processes at arbitrary time instant $t>0$, is obtained. Some numerical results are presented. 
\end{abstract}

\vskip 0.1cm

{\it Keywords:} Telegraph process, telegraph equation, persistent random walk,
probability distribution function of telegraph process,  
Euclidean distance between telegraph processes  

\vskip 0.2cm

{\it AMS 2010 Subject Classification:} 60K35, 60J60, 60J65, 82C41,
82C70

\section{Introduction}

\numberwithin{equation}{section}

The classical telegraph process $X(t)$ is performed by the stochastic motion of a particle that moves
on the real line $\Bbb R$ at some constant finite speed $c$ and alternates two possible directions of motion (forward and backward)
at Poisson-distributed random instants of intensity $\lambda>0$. This random walk was first introduced in the works of Goldstein \cite{gold}
and Kac \cite{kac} (of which the latter is a reprinting of an earlier 1956 article).  The most remarkable fact is that the transition density of
$X(t)$ is the fundamental solution to the hyperbolic telegraph equation (which is one of the classical equations of mathematical physics) and, under
increasing $c$ and $\lambda$, it transforms into the transition density of the standard Brownian motion on $\Bbb R$. Thus, the telegraph process
can be treated as a finite-velocity counterpart of the one-dimensional Brownian motion. The telegraph process $X(t)$ can also be treated in a more
general context of random evolutions (see \cite{pin}).

During last decades the Goldstein-Kac telegraph process and
its numerous generalizations have become the subject of extensive
researches due to great theoretical importance of the model and its 
fruitful applications in statistical physics, financial modeling,
transport phenomena in physical and biological systems, hydrology
and some other fields. Some properties of the solution space of
the Goldstein-Kac telegraph equation were studied in \cite{bart2}. The process of one-dimensional random
motion at finite speed governed by a Poisson process with a
time-dependent parameter was considered in \cite{kap}. The
relationships between the Goldstein-Kac model and physical
processes, including some emerging effects of the relativity
theory, were examined in \cite{bart1}, \cite{cane1}, \cite{cane2}. Formulas for the distributions of the
first-exit time from a given interval and of the maximum
displacement of the telegraph process were obtained in \cite[Section 0.5]{pin}, \cite{foong1}, \cite{mas1}, \cite{mas2}. The behavior of the telegraph process with
absorbing and reflecting barriers was examined in \cite{foong2}, \cite{rat1}. A one-dimensional stochastic
motion with an arbitrary number of velocities and of governing Poisson processes was examined in \cite{kol1}. The
telegraph processes with random velocities were studied in \cite{sta}. The behaviour of the telegraph-type evolutions
in inhomogeneous environments were considered in \cite{rat2}.
Probabilistic methods of solving the Cauchy problems for the telegraph equation
were developed in \cite{kac}, \cite{kis}, \cite{kab}, \cite{turb}. A generalization of the Goldstein-Kac
model for the case of a damped telegraph process with logistic stationary distributions was given in 
\cite{cres2}. A random motion with velocities alternating at Erlang-distributed random times was studied in \cite{cres1}. 
A detailed moment analysis of the telegraph process was done in \cite{kol2}. Explicit formulas for the occupation time distributions of
the telegraph process were recently obtained in \cite{br}.

In all of these works the main subject of interest was a single
particle moving in random direction with constant finite speed. 
In the present article we start studying the problems devoted to
the evolutions of several particles moving randomly with finite
speed. This type of problems can be of a great interest from the
point of view of various possible interactions of the particles.
Such many-particle random motions with interactions can serve as
very good and adequate mathematical models for describing various
real phenomena in physics, chemistry, biology, financial markets
and other fields. For example, in physics and chemistry the
particles are the atoms or molecules of the substance and their
interaction can provoke a physical or chemical reaction. In 
biology the particles can be imagined as the biological objects
(cells, bacteria, animals etc.) and their "interaction" can mean
creating a new cell or, contrary, killing the cell, launching the
infection mechanism or founding a new animal population, respectively. In
financial markets the moving particles can be interpreted as oscillating exchange
rates or stock prices and their "interaction" can mean gaining or ruining.

The system of two particles, whose evolutions do not depend each of other, is the basic model because a many-particle
system can be studied by means of successive growing the
two-particle model. That is why in this article we concentrate our
attention on the basic system of two particles that move, independently each of other, with some constant finite speeds $c_1>0, \; c_2>0,$
on the real line $\Bbb R$ and whose evolutions are driven
by two independent homogeneous Poisson processes of rates $\lambda_1 > 0, \; \lambda_2>0$, respectively.

Let $X_1(t)$ and $X_2(t)$ be two telegraph processes representing the positions of these particles on
$\Bbb R$ at an arbitrary time instant $t>0$. In describing the phenomena of interaction the Euclidean distance 
between these processes  
\begin{equation}\label{intr1}
\rho(t) = \vert X_1(t) - X_2(t) \vert , \qquad t>0.
\end{equation}
is of a special importance. It is quite natural to consider that the particles do not "feel" each 
other if $\rho(t)$ is large. In other words, the forces acting
between the particles are negligible if the distance $\rho(t)$ is
sufficiently big. However, as soon as the distance between the
particles becomes less than some given $r>0$, the particles can start interacting with some positive probability. 
This means that the occurrence of the random event $\{ \rho(t)< r \} $ is the
necessary (but, maybe, not sufficient) condition for launching the
process of interaction at time instant $t>0$. Therefore, the 
distribution $\text{Pr} \{ \rho(t)< r \}$ plays a very important role in analyzing such processes.

The paper is organized as follows. In Section 2 we remind some basic properties of the telegraph process
$X(t)$ that we will substantially be relying on. In Section 3 we obtain a series representation of the probability of being, at time $t>0$,
in an arbitrary subinterval of the support of $X(t)$ and derive a closed-form expression for the probability distribution function
of $X(t)$ that, to the best of our knowledge, were not obtained in the literature. These results are given in terms of
Gauss hypergeometric functions, as well as in terms of Gegenbauer polynomials with non-integer negative upper indices. In Section 4 we formulate 
and prove the principal result of the article yielding the closed-form expression for the probability distribution function of 
the Euclidean distance (\ref{intr1}) between two independent telegraph processes. Some approximate numerical results related to the formula obtained 
are presented in Section 5. In Appendix we prove two auxiliary lemmas related to some indefinite integrals 
of modified Bessel functions and conditional probabilities that are used in our analysis.

\section{Some Basic Properties of Telegraph Process}

\numberwithin{equation}{section}

The telegraph stochastic process is performed by a particle that
starts at the initial time instant $t=0$ from the origin $x=0$ of the real
line $\Bbb R$ and moves with some finite constant speed $c$. The
initial direction of the motion (positive or negative) is taken on
with equal probabilities 1/2. The motion is driven by a
homogeneous Poisson process $N(t)$ of rate $\lambda>0$ as follows. As a
Poisson event occurs, the particle instantaneously takes on the
opposite direction and keeps moving with the same speed $c$ until
the next Poisson event occurrence, then it takes on the opposite
direction again independently of its previous motion, and so on.
This random motion was first studied by Goldstein \cite{gold}
and Kac \cite{kac} and was called the {\it telegraph process} afterwards.

Let $X(t)$ denote the particle's position on $\Bbb R$ at an
arbitrary time instant $t>0$. Since the speed $c$ is finite, then,
at instant $t>0$, the distribution $\text{Pr}\{ X(t)\in
dx \}$ is concentrated in the finite interval $[-ct, ct]$ which is
the support of the distribution of $X(t)$. The density $f(x,t),
\; x\in\Bbb R, \; t\ge 0,$ of the distribution $\text{Pr}\{
X(t)\in dx \}$ has the structure
$$f(x, t) = f^{s}(x, t) + f^{ac}(x, t),$$
where $f^{s}(x, t)$ and $f^{ac}(x, t)$ are the densities of the
singular (with respect to the Lebesgue measure on the line) and of
the absolutely continuous components of the distribution of
$X(t)$, respectively.

The singular component of the distribution is, obviously,
concentrated in two terminal points $\pm ct$ of the interval
$[-ct, ct]$ and corresponds to the case when no one Poisson event
occurs until the moment $t$ and, therefore, the particle does not
change its initial direction. Therefore, the probability of being,
at arbitrary instant $t>0$, at the terminal points $\pm ct$ is
\begin{equation}\label{prop1}
\text{Pr}\left\{ X(t) = ct \right\} = \text{Pr}\left\{ X(t)
= - ct \right\} = \frac{1}{2} \; e^{-\lambda t} .
\end{equation}
The absolutely continuous component of the distribution of $X(t)$
is concentrated in the open interval $(-ct, ct)$ and corresponds
to the case when at least one Poisson event occurs by the moment
$t$ and, therefore, the particle changes its initial direction.
The probability of this event is
\begin{equation}\label{prop2}
\text{Pr}\left\{ X(t) \in (-ct, ct) \right\} = 1 - e^{-\lambda t}.
\end{equation}

The principal result by Goldstein \cite{gold} and Kac \cite{kac}
states that the density $f = f(x,t), \; x\in [-ct, ct], \; t>0,$
of the distribution of $X(t)$ satisfies the following hyperbolic
partial differential equation
\begin{equation}\label{prop3}
\frac{\partial^2 f}{\partial t^2} + 2\lambda
\frac{\partial f}{\partial t} - c^2 \frac{\partial^2 f}{\partial x^2} = 0,
\end{equation}
(which is referred to as the {\it telegraph} or {\it damped wave}
equation) and can be found by solving (\ref{prop3}) with the initial
conditions
\begin{equation}\label{iprop3}
f(x,t)\vert_{t=0} = \delta(x), \qquad
\left.\frac{\partial f(x,t)}{\partial t}\right\vert_{t=0} = 0,
\end{equation}
where $\delta(x)$ is the Dirac delta-function. This means that the
transition density $f(x,t)$ of the process $X(t)$ is the
fundamental solution (i.e. the Green's function) to the telegraph
equation (\ref{prop3}).

The explicit form of the density $f(x,t)$ is given by the formula
(see, for instance, \cite[Section 0.4]{pin}:
\begin{equation}\label{prop4}
\aligned
f(x,t) & = \frac{e^{-\lambda t}}{2} \left[ \delta(ct-x) + \delta(ct+x) \right]\\
& \quad + \frac{\lambda e^{-\lambda t}}{2c} \left[ I_0\left( \frac{\lambda}{c} \sqrt{c^2t^2-x^2} \right) +
\frac{ct}{\sqrt{c^2t^2-x^2}} I_1\left( \frac{\lambda}{c} \sqrt{c^2t^2-x^2} \right) \right] \Theta(ct-\vert x\vert),
\endaligned
\end{equation}
where $I_0(z)$ and $I_1(z)$ are the modified Bessel functions of zero and first orders, respectively, (that is, the
Bessel functions with imaginary argument) given by the formulas
\begin{equation}\label{pprop5}
I_0(z) = \sum_{k=0}^{\infty} \frac{1}{(k!)^2} \left( \frac{z}{2} \right)^{2k} \qquad
I_1(z) = \sum_{k=0}^{\infty} \frac{1}{k! \; (k+1)!} \left( \frac{z}{2} \right)^{2k+1} ,
\end{equation}
and $\Theta(x)$ is the Heaviside step function
\begin{equation}\label{pprop4}
\Theta(x) = \left\{ \aligned 1, \qquad  & \text{if} \; x>0,\\
                               0, \qquad & \text{if} \; x\le 0.
\endaligned \right.
\end{equation}

The first term in (\ref{prop4})
\begin{equation}\label{prop5}
f^s(x,t) = \frac{e^{-\lambda t}}{2} \left[ \delta(ct-x) +
\delta(ct+x) \right]
\end{equation}
is the singular part of the density of the distribution
of $X(t)$ concentrated at two terminal points $\pm ct$ of the
interval $[-ct, ct]$, while the second term in (\ref{prop4})
\begin{equation}\label{prop6}
f^{ac}(x,t) = \frac{\lambda e^{-\lambda t}}{2c}
\left[ I_0\left( \frac{\lambda}{c} \sqrt{c^2t^2-x^2} \right) +
\frac{ct}{\sqrt{c^2t^2-x^2}} I_1\left( \frac{\lambda}{c}
\sqrt{c^2t^2-x^2} \right) \right] \Theta(ct-\vert x\vert),
\end{equation}
represents the density of the absolutely continuous part of the distribution of $X(t)$ concentrated in the open interval $(-ct,ct)$. 

\section{Distribution Function of Telegraph Process}

\numberwithin{equation}{section}

Consider a telegraph process $X(t)$ performed by the stochastic motion of a particle starting, at the initial time instant $t=0$, 
from the origin $x=0$ of the real line $\Bbb R$ and moving with some finite constant speed $c>0$ whose evolution is driven 
by a homogeneous Poisson process $N(t)$ of rate $\lambda >0$, as described above. 

As is noted above, at arbitrary time instant $t>0$ the process $X(t)$ is concentrated in the interval $[-ct, ct]$. 
Let $a, b \in \Bbb R, \; a < b$, be arbitrary points of $\Bbb R$ such that the intervals $(a, b)$ and $(-ct, ct)$ have a non-empty intersection, that is, $(a, b)\cap (-ct, ct) \neq \varnothing$. We are interested in the probability $\text{Pr} \left\{ X(t)\in (a, b)\cap (-ct, ct) \right\}$ that the process $X(t)$, at time instant $t>0$, is located in the subinterval $(a, b)\cap (-ct, ct)\subseteq (-ct,ct)$. This result is presented by the following proposition.

\bigskip

{\bf Proposition 1.} {\it For arbitrary time instant $t>0$ and arbitrary open interval $(a, b)\subset\Bbb R$, $a, b \in \Bbb R, \; a < b$, such that $(a, b)\cap (-ct, ct) \neq \varnothing$ the following formula holds:}
\begin{equation}\label{sub1}
\aligned
& \text{Pr} \left\{ X(t)\in (a, b)\cap (-ct, ct) \right\} \\ 
& = \frac{\lambda e^{-\lambda t}}{2c} \sum_{k=0}^{\infty} \frac{1}{(k!)^2} \left(\frac{\lambda t}{2}\right)^{2k} \left( 1+ \frac{\lambda t}{2k+2} \right) 
\left[ \beta F\left( -k, \frac{1}{2}; \frac{3}{2}; \frac{\beta^2}{c^2t^2} \right) - \alpha F\left( -k, \frac{1}{2}; \frac{3}{2}; \frac{\alpha^2}{c^2t^2} \right) \right] , 
\endaligned
\end{equation}
{\it where}
\begin{equation}\label{subb1}
\alpha = \max\{ -ct, \; a\}, \qquad \beta = \min\{ ct, \; b \} 
\end{equation}
{\it and}
$$F(\xi,\eta; \zeta; z) = \; _{2}F_1(\xi,\eta; \zeta; z) = \sum_{k=0}^{\infty} \frac{(\xi)_k \; (\eta)_k}{(\zeta)_k} \; \frac{z^k}{k!}$$
{\it is the Gauss hypergeometric function.}

\vskip 0.1cm

{\it Proof.} By integrating density (\ref{prop6}) and applying formulas (\ref{app4}) and (\ref{app5}) of the Appendix (see below), we obtain:
$$\aligned
\text{Pr} & \bigl\{ X(t)\in (a, b)\cap (-ct, ct) \bigr\} \\ 
& = \frac{\lambda e^{-\lambda t}}{2c} \left[ \; \int\limits_{\alpha}^{\beta} I_0\left( \frac{\lambda}{c} \sqrt{c^2t^2-x^2} \right) \; dx
+ ct \int\limits_{\alpha}^{\beta} \frac{I_1\left( \frac{\lambda}{c}
\sqrt{c^2t^2-x^2} \right)}{\sqrt{c^2t^2-x^2}} \; dx \right] \\
& = \frac{\lambda e^{-\lambda t}}{2c} \left[ x \sum_{k=0}^{\infty} \frac{1}{(k!)^2} \; \left(\frac{\lambda t}{2}\right)^{2k} F\left( -k, \frac{1}{2}; \frac{3}{2};  \frac{x^2}{c^2t^2} \right) \right. \\
& \left.\left. \qquad\qquad + x \sum_{k=0}^{\infty} \frac{1}{k! \; (k+1)!} \; \left(\frac{\lambda t}{2}\right)^{2k+1} F\left( -k, \frac{1}{2}; \frac{3}{2};  \frac{x^2}{c^2t^2} \right)  \right] \right|_{x=\alpha}^{x=\beta} \\
& = \frac{\lambda e^{-\lambda t}}{2c} \sum_{k=0}^{\infty} \frac{1}{(k!)^2} \; \left(\frac{\lambda t}{2}\right)^{2k} \left( 1+ \frac{\lambda t}{2k+2} \right) \\
& \qquad\qquad\qquad \times \left[ \beta F\left( -k, \frac{1}{2}; \frac{3}{2}; \frac{\beta^2}{c^2t^2} \right) - \alpha F\left( -k, \frac{1}{2}; \frac{3}{2};  \frac{\alpha^2}{c^2t^2} \right) \right] ,
\endaligned$$
proving (\ref{sub1}). $\square$

\bigskip

{\it Remark 1.}  Let $x\in (-ct, ct)$ be an arbitrary interior point of the open interval $(-ct, ct)$ and let $r>0$ be an arbitrary 
positive number such that $(x-r, x+r) \subseteq (-ct, ct)$. Then, according to (\ref{sub1}), we obtain the following formula for 
the probability of being in the subinterval $(x-r, x+r)\subseteq (-ct, ct)$ of radius $r$ centered at point $x$:
\begin{equation}\label{ssubb2}
\aligned 
\text{Pr} & \left\{ X(t)\in (x-r, x+r) \right\} \\
& = \frac{\lambda e^{-\lambda t}}{2c} \sum_{k=0}^{\infty} \frac{1}{(k!)^2} \; \left( \frac{\lambda t}{2} \right)^{2k} 
\left( 1 + \frac{\lambda t}{2k+2} \right) \\
& \qquad\qquad \times \left[ (x+r) F\left( -k, \frac{1}{2}; \frac{3}{2}; \frac{(x+r)^2}{c^2t^2} \right) - (x-r) F\left( -k, \frac{1}{2}; \frac{3}{2}; \frac{(x-r)^2}{c^2t^2} \right) \right],
\endaligned 
\end{equation}
$$-ct \le x-r < x+r \le ct .$$
By setting $x=0$ in (\ref{ssubb2}) we obtain the formula: 
\begin{equation}\label{subb2}
\text{Pr} \left\{ X(t)\in (-r, r) \right\} = \frac{\lambda r \; e^{-\lambda t}}{c} \sum_{k=0}^{\infty} \frac{1}{(k!)^2} \left( \frac{\lambda t}{2} \right)^{2k} 
\left( 1 + \frac{\lambda t}{2k+2} \right) F\left( -k, \frac{1}{2}; \frac{3}{2}; \frac{r^2}{c^2t^2} \right) ,
\end{equation}
yielding the probability of being in the symmetric (with respect to the start point $x=0$) subinterval $(-r, r)\subseteq (-ct, ct)$. In particular, setting $r=ct$ in (\ref{subb2}) and applying formula (\ref{gauss1}) (see below) we get 
$$\aligned
\text{Pr} \left\{ X(t)\in (-ct, ct) \right\} & = \lambda t \; e^{-\lambda t} \sum_{k=0}^{\infty} \frac{1}{(k!)^2} \left( \frac{\lambda t}{2} \right)^{2k}
\left( 1 + \frac{\lambda t}{2k+2} \right) F\left( -k, \frac{1}{2}; \frac{3}{2}; 1 \right) \\
& = \lambda t \; e^{-\lambda t} \sum_{k=0}^{\infty} \frac{1}{(k!)^2} \left( \frac{\lambda t}{2} \right)^{2k}
\left( 1 + \frac{\lambda t}{2k+2} \right) \; \frac{2^k \; k!}{(2k+1)!!} \\
& = e^{-\lambda t} \sum_{k=0}^{\infty} \frac{(\lambda t)^{2k+1}}{2^k \; k! \; (2k+1)!!} \left( 1 + \frac{\lambda t}{2k+2} \right) \\
& = e^{-\lambda t} \sum_{k=0}^{\infty} \frac{(\lambda t)^{2k+1}}{(2k+1)!} \left( 1 + \frac{\lambda t}{2k+2} \right) \\
& = e^{-\lambda t} \left[ \sum_{k=0}^{\infty} \frac{(\lambda t)^{2k+1}}{(2k+1)!} + \sum_{k=0}^{\infty} \frac{(\lambda t)^{2k+2}}{(2k+2)!} \right] \\
& = e^{-\lambda t} \left[ \sinh(\lambda t) + \cosh(\lambda t) - 1 \right] \\ 
& = e^{-\lambda t} \left[ e^{\lambda t} - 1 \right] \\ 
& = 1 - e^{-\lambda t} ,
\endaligned$$
exactly coinciding with (\ref{prop2}). 

From Proposition 1 we can extract the explicit form of the probability distribution function of $X(t)$.

\bigskip

{\bf Proposition 2.} {\it The probability distribution function of the telegraph process $X(t)$ has the form:} 

\begin{equation}\label{sub2}
\aligned
\text{Pr} & \left\{ X(t) < x \right\} \\
& = \left\{ \aligned
& 0, \qquad\qquad\qquad\qquad\qquad\qquad x\in (-\infty, -ct],\\
& \frac{1}{2} + \frac{\lambda x  e^{-\lambda t}}{2c} \sum_{k=0}^{\infty} \frac{1}{(k!)^2} \left(\frac{\lambda t}{2}\right)^{2k} \left( 1+ \frac{\lambda t}{2k+2} \right) F\left( -k, \frac{1}{2}; \frac{3}{2}; \frac{x^2}{c^2t^2} \right) , \qquad  x\in (-ct, ct],\\
& 1, \qquad\qquad\qquad\qquad\qquad\qquad x\in (ct, +\infty).\\
\endaligned \right. \endaligned
\end{equation}

\vskip 0.1cm

{\it Proof.} According to Proposition 1, for arbitrary $x\in (-ct, ct)$ we have
$$\aligned
& \text{Pr} \left\{ X(t) \in (-ct, x) \right\} \\ 
& = \frac{\lambda e^{-\lambda t}}{2c} \sum_{k=0}^{\infty} \frac{1}{(k!)^2} \left(\frac{\lambda t}{2}\right)^{2k} \left( 1 + \frac{\lambda t}{2k+2} \right) \left[ x F\left( -k, \frac{1}{2}; \frac{3}{2}; \frac{x^2}{c^2t^2} \right) + ct F\left( -k, \frac{1}{2}; \frac{3}{2}; 1  \right) \right] \\
& = \frac{\lambda x e^{-\lambda t}}{2c} \sum_{k=0}^{\infty} \frac{1}{(k!)^2} \left(\frac{\lambda t}{2}\right)^{2k} \left( 1 + \frac{\lambda t}{2k+2} \right) F\left( -k, \frac{1}{2}; \frac{3}{2}; \frac{x^2}{c^2t^2} \right) \\ 
& \qquad\qquad + e^{-\lambda t} \sum_{k=0}^{\infty} \frac{1}{(k!)^2} \left(\frac{\lambda t}{2}\right)^{2k+1} \left( 1 + \frac{\lambda t}{2k+2} \right) F\left( -k, \frac{1}{2}; \frac{3}{2}; 1 \right) .
\endaligned$$
Consider separately the second term of this expression. Taking into account that (see also general formula (\ref{gauss2}) below)
\begin{equation}\label{gauss1}
F\left( -k, \frac{1}{2}; \frac{3}{2}; 1 \right) = \frac{(2k)!!}{(2k+1)!!} = \frac{2^k \; k!}{(2k+1)!!} , \qquad k\ge 0, 
\end{equation}
we get  
$$\aligned
e^{-\lambda t} \sum_{k=0}^{\infty} & \frac{1}{(k!)^2} \left(\frac{\lambda t}{2}\right)^{2k+1} \left( 1 + \frac{\lambda t}{2k+2} \right) F\left( -k, \frac{1}{2}; \frac{3}{2}; 1  \right) \\ 
& = e^{-\lambda t} \sum_{k=0}^{\infty} \frac{2^k}{k! \; (2k+1)!!} \left(\frac{\lambda t}{2}\right)^{2k+1} \left( 1 + \frac{\lambda t}{2k+2} \right) \\
& = e^{-\lambda t} \sum_{k=0}^{\infty} \frac{2^{2k}}{2^k k! \; (2k+1)!!} \left(\frac{\lambda t}{2}\right)^{2k+1} \left( 1 + \frac{\lambda t}{2k+2} \right) \\
& = \frac{e^{-\lambda t}}{2} \sum_{k=0}^{\infty} \frac{(\lambda t)^{2k+1}}{(2k+1)!} \left( 1 + \frac{\lambda t}{2k+2} \right) \\
& = \frac{e^{-\lambda t}}{2} \left[ \sum_{k=0}^{\infty} \frac{(\lambda t)^{2k+1}}{(2k+1)!} + \sum_{k=0}^{\infty} \frac{(\lambda t)^{2k+2}}{(2k+2)!}  \right] \\
& = \frac{e^{-\lambda t}}{2} \left[ \sinh{(\lambda t)} + \cosh{(\lambda t)} - 1 \right] \\
& = \frac{1}{2} - \frac{e^{-\lambda t}}{2} .
\endaligned$$
Therefore, for arbitrary $x\in (-ct, ct]$ we obtain 

$$\aligned
& \text{Pr} \left\{ X(t) < x \right\} \\
& = \text{Pr} \left\{ X(t) = - ct \right\} + \text{Pr} \left\{ X(t) \in (-ct, x) \right\} \\
& = \frac{e^{-\lambda t}}{2} + \frac{\lambda x e^{-\lambda t}}{2c} \sum_{k=0}^{\infty} \frac{1}{(k!)^2} \left(\frac{\lambda t}{2}\right)^{2k} \left( 1 + \frac{\lambda t}{2k+2} \right) F\left( -k, \frac{1}{2}; \frac{3}{2}; \frac{x^2}{c^2t^2} \right) + \frac{1}{2} - \frac{e^{-\lambda t}}{2} \\
& = \frac{1}{2} + \frac{\lambda x e^{-\lambda t}}{2c} \sum_{k=0}^{\infty} \frac{1}{(k!)^2} \left(\frac{\lambda t}{2}\right)^{2k} \left( 1 + \frac{\lambda t}{2k+2} \right) F\left( -k, \frac{1}{2}; \frac{3}{2}; \frac{x^2}{c^2t^2} \right) .
\endaligned$$
The proposition is proved. $\square$ 

\bigskip

The shape of probability distribution function (\ref{sub2}) in the interval $(-2,2]$ is presented in Figure 1.

\begin{center}
\begin{figure}[htbp]
\centerline{\includegraphics[width=10cm,height=8cm]{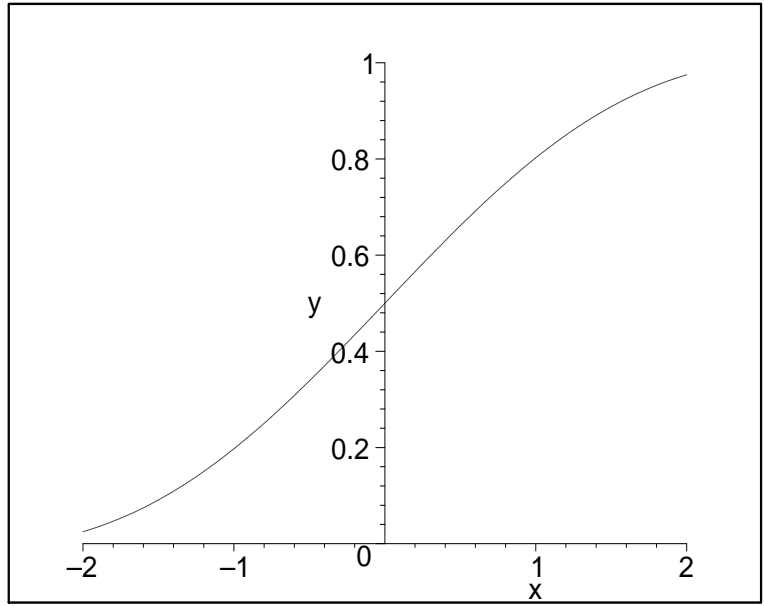}}
\caption{\it The shape of p.d.f. (\ref{sub2}) at instant $t=2$ (for $c=1, \; \lambda =1.5$)}
\end{figure}
\end{center}

{\it Remark 2.} Note that, in view of \cite[page 465, Formula 163]{pbm}, we have
\begin{equation}\label{gauss2}
F\left( -k, \frac{1}{2}; \frac{3}{2}; z \right) = - \frac{(2k)!!}{(2k+1)!! \; \sqrt{z}} \; C_{2k+1}^{-k-1/2} (\sqrt{z}) ,
\end{equation}
where $C_n^{\nu}(z)$ are the Gegenbauer polynomials. Therefore, formulas (\ref{sub1}) and (\ref{sub2}) can be represented in the 
following alternative forms:
\begin{equation}\label{sub4}
\aligned
\text{Pr} & \left\{ X(t)\in (a, b)\cap (-ct, ct) \right\} \\
& = \frac{e^{-\lambda t}}{2} \sum_{k=0}^{\infty} \frac{(\lambda
t)^{2k+1}}{(2k+1)!} \left( 1 + \frac{\lambda t}{2k+2} \right)
\left[ \text{sgn}(\alpha) \; C_{2k+1}^{-k-1/2} \left( \frac{|\alpha|}{ct} \right) - \text{sgn}(\beta) \;
C_{2k+1}^{-k-1/2} \left( \frac{|\beta|}{ct} \right) \right] ,
\endaligned
\end{equation}
where $\alpha$ and $\beta$ are given by (\ref{subb1}) and, respectively, 
\begin{equation}\label{sub5}
\aligned
\text{Pr} & \left\{ X(t) < x \right\} \\
& = \left\{ \aligned
0, \qquad\qquad\qquad\qquad\qquad & x\in (-\infty, -ct],\\
\frac{1}{2} - \frac{e^{-\lambda t}}{2} \; \text{sgn}(x) \sum_{k=0}^{\infty}
\frac{(\lambda t)^{2k+1}}{(2k+1)!} \left( 1 + \frac{\lambda
t}{2k+2} \right) C_{2k+1}^{-k-1/2} \left(
\frac{|x|}{ct} \right), \qquad  & x\in (-ct, ct],\\
1, \qquad\qquad\qquad\qquad\qquad & x\in (ct, +\infty) .
\endaligned \right. \endaligned
\end{equation}

\bigskip

{\it Remark 3.} We see that function (\ref{sub2}) has discontinuities at the points $\pm ct$ determined by the singularities concentrated 
at these two points. It is easy to check that distribution function (\ref{sub2}) produces the expected equalities:
$$\lim\limits_{\varepsilon\to 0+0}  \text{Pr} \left\{ X(t) < -ct + \varepsilon \right\} = \frac{e^{-\lambda t}}{2} , \qquad  
\text{Pr} \left\{ X(t) < ct \right\} = 1 -\frac{e^{-\lambda t}}{2} .$$
This means that probability distribution function (\ref{sub2}) is left-continuous and it has jumps at the terminal points $\pm ct$
of the same amplitude $e^{-\lambda t}/2$.

\section{Euclidean Distance Between Two Telegraph Processes}

\numberwithin{equation}{section}

Consider two independent telegraph processes $X_1(t)$ and $X_2(t)$ performed by the stochastic motions of two particles (as described in Section 2 above) 
with finite speeds $c_1>0, \; c_2>0$ and driven by two independent Poisson processes $N_1(t)$ and $N_2(t)$ of rates $\lambda_1>0, \; \lambda_2>0$, 
respectively. Without loss of generality, we can suppose that, at the initial time instant $t=0$, both the processes $X_1(t)$ and $X_2(t)$ 
simultaneously start from the origin $x=0$ of the real line $\Bbb R$. For the sake of definiteness, we also suppose that $c_1\ge c_2$ 
(otherwise, one can merely change numeration of the processes).   

The subject of our interest is the Euclidean distance 
\begin{equation}\label{dist1}
\rho(t) =  \left| X_1(t) - X_2(t) \right| , \qquad t>0,
\end{equation}
between these processes at time instant $t>0$. 

It is clear that $0 \le\rho(t)\le (c_1+c_2)t$, that is, the interval $[0, (c_1+c_2)t]$ is the support of the distribution $\text{Pr} \left\{ \rho(t) < r \right\}$ 
of process (\ref{dist1}). The distribution of $\rho(t), \; t>0,$ consists of two components. The singular part of the distribution is concentrated at two points 
$(c_1-c_2)t$ and $(c_1+c_2)t$ of the support. For arbitrary $t>0$ the process $\rho(t)$ is located at point $(c_1-c_2)t$ if and only if both the particles take 
the same initial direction (the probability of this event is $1/2$) and no one Poisson event occurs till time instant $t$ (the probability of this event is  $e^{-(\lambda_1+\lambda_2)t}$). Similarly, $\rho(t)$ is located at point $(c_1+c_2)t$ if and only if the particles take different  
initial directions (the probability of this event is $1/2$) and no one Poisson event occurs till time instant $t$ (the probability of this event is  $e^{-(\lambda_1+\lambda_2)t}$). Thus, we have:
\begin{equation}\label{dist2}
\aligned
\text{Pr} \left\{ \rho(t) = (c_1-c_2)t  \right\} & = \frac{1}{2} \; e^{-(\lambda_1+\lambda_2)t} , \\
\text{Pr} \left\{ \rho(t) = (c_1+c_2)t  \right\} & = \frac{1}{2} \; e^{-(\lambda_1+\lambda_2)t} ,
\endaligned \qquad t>0.
\end{equation}
Therefore, the singular part $\varphi^s(r,t)$ of the density $\varphi(r,t)$ of the distribution $\text{Pr} \left\{ \rho(t) < r \right\}$ is the generalized function:
\begin{equation}\label{dist3}
\varphi^s(r,t) = \frac{e^{-(\lambda_1+\lambda_2)t}}{2} \; \bigl[ \delta(r-(c_1-c_2)t) + \delta(r-(c_1+c_2)t) \bigr] , \qquad r\in\Bbb R, \quad t>0, 
\end{equation}
where $\delta(x)$ is the Dirac delta-function. 

The remaining part of distribution is concentrated in the area 
$$M_t = (0, (c_1-c_2)t)\cup ((c_1-c_2)t, (c_1+c_2)t), \qquad t>0,$$ 
(note that if $c_1=c_2=c$ then $M_t$ transforms into the interval $(0, 2ct)$). This is the support of the absolutely continuous part of the distribution 
$\text{Pr} \left\{ \rho(t) < r \right\}$ corresponding to the case when at least one Poisson event occurs till time instant $t>0$. 

Our goal is to obtain an explicit formula for the probability distribution function  
\begin{equation}\label{distPhi}
\Phi(r,t) = \text{Pr} \left\{ \rho(t) < r \right\}, \qquad r\in\Bbb R , \quad t>0, 
\end{equation} 
of the Euclidean distance $\rho(t)$. The form of this distribution function is somewhat different for the cases $c_1=c_2$ and $c_1>c_2$ due to the fact that, if $c_1=c_2$ 
then the singularity point $(c_1-c_2)t=0$ and this is the {\it terminal} point, while in the case $c_1>c_2$ this point is an {\it interior} point of the support. That is why in the following theorem we derive the probability distribution function in the more complicated case $c_1>c_2$. Similar results concerning the more simple case $c_1=c_2$ will be given separately at the end of this section. 

\bigskip

{\bf Theorem 1.} {\it Under the condition $c_1>c_2$ the probability distribution function} (\ref{distPhi}) {\it has the form:}
\begin{equation}\label{dist4}
\Phi(r,t) = \left\{ \aligned 0 , \qquad & \text{if} \; r\in (-\infty, \; 0] ,\\ 
                             G(r,t), \qquad & \text{if} \; r\in (0, \; (c_1-c_2)t] ,\\
                             Q(r,t), \qquad & \text{if} \; r\in ((c_1-c_2)t, \; (c_1+c_2)t] ,\\
                             1 , \qquad & \text{if} \; r\in ((c_1+c_2)t, \; +\infty) ,
\endaligned \right.
\end{equation}
$$\qquad r\in\Bbb R, \qquad t>0, \qquad c_1>c_2,$$
{\it where functions $G(r,t)$ and $Q(r,t)$ are given by the formulas:}
\begin{equation}\label{distG}
\aligned 
G(r,t) & = \frac{\lambda_1 \; e^{-(\lambda_1+\lambda_2)t}}{2c_1} \sum_{k=0}^{\infty} \frac{1}{(k!)^2} \; \left( \frac{\lambda_1 t}{2} \right)^{2k} \left( 1 + \frac{\lambda_1 t}{2k+2} \right) \\ 
& \qquad \times \biggl[ (c_2t+r) F\left( -k, \frac{1}{2}; \; \frac{3}{2}; \; \frac{(c_2t+r)^2}{c_1^2t^2} \right) - (c_2t-r) F\left( -k, \frac{1}{2}; \; \frac{3}{2}; \; \frac{(c_2t-r)^2}{c_1^2t^2} \right) \biggr] \\ 
& + \frac{\lambda_1 \; e^{-\lambda_1t}}{2c_1} \sum_{k=0}^{\infty} \frac{1}{(k!)^2} \; \left( \frac{\lambda_1 t}{2} \right)^{2k} \left( 1 + \frac{\lambda_1 t}{2k+2} \right)  \mathcal I_k(r,t) , 
\endaligned
\end{equation}

\begin{equation}\label{distQ}
\aligned 
Q(r,t) & = \frac{1}{2} \biggl[ \left( 1-e^{-\lambda_1t} \right) e^{-\lambda_2t} + \left( 1-e^{-\lambda_2t} \right) e^{-\lambda_1t} + e^{-(\lambda_1+\lambda_2)t} \biggr] \\ 
& \; - \frac{\lambda_1 (c_2t-r) e^{-(\lambda_1+\lambda_2)t}}{2c_1} \\ 
& \qquad\qquad\qquad \times \sum_{k=0}^{\infty} \frac{1}{(k!)^2} \; \left( \frac{\lambda_1 t}{2} \right)^{2k} \left( 1 + \frac{\lambda_1 t}{2k+2} \right) F\left( -k, \frac{1}{2}; \; \frac{3}{2}; \; \frac{(c_2t-r)^2}{c_1^2t^2} \right) \\ 
& \; - \frac{\lambda_2 (c_1t-r) e^{-(\lambda_1+\lambda_2)t}}{2c_2} \\ 
& \qquad\qquad\qquad \times \sum_{k=0}^{\infty} \frac{1}{(k!)^2} \; \left( \frac{\lambda_2 t}{2} \right)^{2k} \left( 1 + \frac{\lambda_2 t}{2k+2} \right) F\left( -k, \frac{1}{2}; \; \frac{3}{2}; \; \frac{(c_1t-r)^2}{c_2^2t^2} \right) \\ 
& \; + \frac{\lambda_1 \; e^{-\lambda_1t}}{2c_1} \sum_{k=0}^{\infty} \frac{1}{(k!)^2} \; \left( \frac{\lambda_1 t}{2} \right)^{2k} \left( 1 + \frac{\lambda_1 t}{2k+2} \right) \mathcal I_k(r,t) , 
\endaligned
\end{equation}
{\it with the integral term} 
\begin{equation}\label{distINT}
\aligned 
\mathcal I_k(r,t) & = \frac{\lambda_2 e^{-\lambda_2 t}}{2c_2} \int\limits_{-c_2t}^{c_2t} \biggl[ \beta(x,r) F\left( -k, \frac{1}{2}; \frac{3}{2}; \frac{(\beta(x,r))^2}{c_1^2t^2} \right) - \alpha(x,r) F\left( -k, \frac{1}{2}; \frac{3}{2}; \frac{(\alpha(x,r))^2}{c_1^2t^2} \right) \biggr] \\ 
& \hskip 3cm \times \biggl[ I_0\left( \frac{\lambda_2}{c_2} \sqrt{c_2^2t^2-x^2} \right) + \frac{c_2t}{\sqrt{c_2^2t^2-x^2}} I_1\left( \frac{\lambda_2}{c_2} 
\sqrt{c_2^2t^2-x^2} \right) \biggr] dx ,
\endaligned
\end{equation}
{\it where the variables $\alpha(x,r)$ and $\beta(x,r)$ are defined by the formulas}: 
\begin{equation}\label{distAB}
\alpha(x,r) = \max\{ -c_1t, \; x-r \} , \qquad \beta(x,r) = \min\{ c_1t, \; x+r \} , 
\end{equation}
$$x\in(-c_2t, \; c_2t), \qquad r\in M_t.$$

\vskip 0.2cm

{\it Proof.} For distribution function (\ref{distPhi}) we have 
\begin{equation}\label{dist5}
\aligned
\Phi(r,t) & = e^{-(\lambda_1+\lambda_2)t} \; \text{Pr} \left\{ \rho(t) < r \; \bigl| \; N_1(t)=0, \; N_2(t)=0 \right\} \\
& \quad + \left( 1-e^{-\lambda_1t} \right) \; e^{-\lambda_2t} \; \text{Pr} \left\{ \rho(t) < r \; \bigl| \; N_1(t)\ge 1, \; N_2(t)=0 \right\} \\
& \quad + e^{-\lambda_1t} \; \left( 1- e^{-\lambda_2t} \right) \; \text{Pr} \left\{ \rho(t) < r \; \bigl| \; N_1(t)=0, \; N_2(t)\ge 1 \right\} \\
& \quad + \left( 1 - e^{-\lambda_1t} \right) \; \left( 1- e^{-\lambda_2t} \right) \; \text{Pr} \left\{ \rho(t) < r \; \bigl| \; N_1(t)\ge 1, \; N_2(t)\ge 1 \right\} .
\endaligned
\end{equation}
Let us evaluate separately conditional probabilities on the right-hand side of (\ref{dist5}). Obviously, the first conditional probability is: 
\begin{equation}\label{ddist5}
\text{Pr} \left\{ \rho(t) < r \; \bigl| \; N_1(t)=0, \; N_2(t)=0 \right\} = \left\{ 
\aligned 0, \qquad & \text{if} \; r\in(-\infty, \; (c_1-c_2)t] ,\\
         \frac{1}{2}, \qquad & \text{if} \; r\in((c_1-c_2)t, \; (c_1+c_2)t] ,\\
         1, \qquad & \text{if} \; r\in((c_1+c_2)t, \; +\infty) .
\endaligned \right.
\end{equation}

\vskip 0.1cm

$\bullet$ {\it Evaluation of} $\text{Pr} \left\{ \rho(t) < r \; \bigl| \; N_1(t)\ge 1, \; N_2(t)=0 \right\}$. We note that the following equalities for random events hold:
$$\aligned 
\left\{ N_1(t)\ge 1 \right\} & = \left\{ X_1(t) \in (-c_1t, \; c_1t) \right\} , \\
\left\{ N_2(t)=0 \right\} & = \left\{ X_2(t)=-c_2t \right\} + \left\{ X_2(t)=c_2t \right\} .
\endaligned$$
Therefore, according to formula (\ref{app8}) of Lemma A2, we have 
$$\aligned 
& \text{Pr} \left\{ \rho(t) < r \; \bigl| \; N_1(t)\ge 1, \; N_2(t)=0 \right\} \\
& = \text{Pr} \bigl\{ \rho(t) < r \; \bigl| \; \left\{ X_1(t) \in (-c_1t, \; c_1t) \right\} \cap \left( \left\{ X_2(t)=-c_2t \right\} + \left\{ X_2(t)=c_2t \right\} \right) \bigr\} \\
& = \frac{1}{2} \biggl[  \text{Pr} \bigl\{ \rho(t) < r \; \bigl| \; \left\{ X_1(t) \in (-c_1t, \; c_1t) \right\} \cap \left\{ X_2(t)=-c_2t \right\} \bigr\} \\
& \hskip 2cm + \text{Pr} \bigl\{ \rho(t) < r \; \bigl| \; \left\{ X_1(t) \in (-c_1t, \; c_1t) \right\} \cap \left\{ X_2(t)=c_2t \right\} \bigr\}  \biggr] \\
& = \frac{1}{2} \biggl[  \frac{\text{Pr} \bigl\{ \left\{ X_1(t) \in (X_2(t)-r, \; X_2(t)+r) \right\} \cap  \left\{ X_1(t) \in (-c_1t, \; c_1t) \right\} \cap \left\{ X_2(t)=-c_2t \right\} \bigr\}}{ {\text{Pr}} \left\{ X_1(t) \in (-c_1t, \; c_1t) \right\} \; {\text{Pr}} \left\{ X_2(t) = -c_1t \right\} } \\
& \hskip 1cm + \frac{\text{Pr} \bigl\{ \left\{ X_1(t) \in (X_2(t)-r, \; X_2(t)+r) \right\} \cap  \left\{ X_1(t) \in (-c_1t, \; c_1t) \right\} \cap \left\{ X_2(t)=c_2t \right\} \bigr\}}{ {\text{Pr}} \left\{ X_1(t) \in (-c_1t, \; c_1t) \right\} \; {\text{Pr}} \left\{ X_2(t) = c_1t \right\} }  \biggr] \\
& = \frac{1}{2(1-e^{-\lambda_1t})} \biggl[  \text{Pr} \bigl\{ \left\{ X_1(t) \in (-c_2t-r, \; -c_2t+r) \right\} \cap  \left\{ X_1(t) \in (-c_1t, \; c_1t) \right\} \bigr\} \\
& \hskip 4cm + \text{Pr} \bigl\{ \left\{ X_1(t) \in (c_2t-r, \; c_2t+r) \right\} \cap  \left\{ X_1(t) \in (-c_1t, \; c_1t) \right\} \bigr\} \biggr] \\
& = \frac{1}{2(1-e^{-\lambda_1t})} \biggl[  \text{Pr} \left\{ X_1(t) \in (\alpha, \; -c_2t+r) \right\} + \text{Pr} \left\{ X_1(t) \in (c_2t-r, \; \beta) \right\} \biggr] ,
\endaligned$$
where   
$$\alpha = \max\{-c_1t, \; -c_2t-r \},  \qquad \beta = \min\{c_1t, \; c_2t+r \} .$$
Applying formula (\ref{sub1}) of Proposition 1 we get:
\begin{equation}\label{dist6}
\aligned
& \text{Pr} \left\{ \rho(t) < r \; \bigl| \; N_1(t)\ge 1, \; N_2(t)=0 \right\} \\
& = \frac{1}{2(1-e^{-\lambda_1t})} \biggl\{ \frac{\lambda_1 e^{-\lambda_1 t}}{2c_1} \sum_{k=0}^{\infty} \frac{1}{(k!)^2} \; \left(\frac{\lambda_1 t}{2}\right)^{2k} \left( 1 + \frac{\lambda_1 t}{2k+2} \right) \\
& \qquad\qquad\qquad \times \biggl[ \beta  F\left( -k, \frac{1}{2}; \; \frac{3}{2}; \; \frac{\beta^2}{c_1^2t^2} \right) - (c_2t-r) F\left( -k, \frac{1}{2}; \; \frac{3}{2}; \; \frac{(c_2t-r)^2}{c_1^2t^2} \right) \biggr] \\
& \qquad\qquad\qquad + \frac{\lambda_1 e^{-\lambda_1 t}}{2c_1} \sum_{k=0}^{\infty} \frac{1}{(k!)^2} \; \left(\frac{\lambda_1 t}{2}\right)^{2k} \left( 1+ \frac{\lambda_1 t}{2k+2} \right) \\
& \qquad\qquad\qquad \times \biggl[ (-c_2t+r) F\left( -k, \frac{1}{2}; \; \frac{3}{2}; \; \frac{(-c_2t+r)^2}{c_1^2t^2} \right) - \alpha F\left( -k, \frac{1}{2}; \; \frac{3}{2}; \; \frac{\alpha^2}{c_1^2t^2} \right) \biggr] \biggr\} .
\endaligned
\end{equation}
It is easy to check that 
$$\beta = \left\{ \aligned c_2t+r, \qquad & \text{if} \; r\in (0, \; (c_1-c_2)t] ,\\ 
                                          c_1t,   \qquad & \text{if} \; r\in ((c_1-c_2)t, \; (c_1+c_2)t] ,  
\endaligned \right. $$
$$\alpha = \left\{ \aligned -c_2t-r, \qquad & \text{if} \; r\in (0, \; (c_1-c_2)t] ,\\ 
                                             -c_1t,   \qquad & \text{if} \; r\in ((c_1-c_2)t, \; (c_1+c_2)t] .                                
\endaligned \right. $$
From these formulas we see that $\alpha = -\beta$ independently of $r$. Therefore, (\ref{dist6}) becomes: 
\begin{equation}\label{ddist6}
\aligned
\text{Pr} & \bigl\{ \rho(t) < r \; \bigl| \; N_1(t)\ge 1, \; N_2(t)=0 \bigr\} \\
& = \frac{\lambda_1 \; e^{-\lambda_1 t}}{2c_1 (1-e^{-\lambda_1t})} \sum_{k=0}^{\infty} \frac{1}{(k!)^2} \; \left(\frac{\lambda_1 t}{2}\right)^{2k} \left( 1+ \frac{\lambda_1 t}{2k+2} \right) \\
& \qquad\qquad\qquad \times \biggl[ \beta  F\left( -k, \frac{1}{2}; \; \frac{3}{2}; \; \frac{\beta^2}{c_1^2t^2} \right) - (c_2t-r) F\left( -k, \frac{1}{2}; \; \frac{3}{2}; \; \frac{(c_2t-r)^2}{c_1^2t^2} \right) \biggr] .
\endaligned
\end{equation}

If $r\in (0, \; (c_1-c_2)t]$ then $\beta=c_2t+r$ and, therefore formula (\ref{ddist6}) takes the form: 
\begin{equation}\label{dist7}
\aligned
\text{Pr} & \bigl\{ \rho(t) < r \; \bigl| \; N_1(t)\ge 1, \; N_2(t) = 0 \bigr\} \\
& = \frac{\lambda_1 \; e^{-\lambda_1 t}}{2c_1(1-e^{-\lambda_1t})} \sum_{k=0}^{\infty} \frac{1}{(k!)^2} \; \left(\frac{\lambda_1 t}{2}\right)^{2k} \left( 1+ \frac{\lambda_1 t}{2k+2} \right) \\
& \qquad \times \left[ (c_2t+r) \; F\left( -k, \frac{1}{2}; \; \frac{3}{2}; \; \frac{(c_2t+r)^2}{c_1^2t^2} \right) - (c_2t-r) \; F\left( -k, \frac{1}{2}; \; \frac{3}{2}; \; \frac{(c_2t-r)^2}{c_1^2t^2} \right)  \right] , 
\endaligned
\end{equation}
$$\text{if} \; r\in (0, \; (c_1-c_2)t] .$$

For $r\in ((c_1-c_2)t, \; (c_1+c_2)t]$ formula (\ref{ddist6}) gives:  
\begin{equation}\label{dist8}
\aligned
\text{Pr} & \bigl\{ \rho(t) < r \; \bigl| \; N_1(t)\ge 1, \; N_2(t) = 0 \bigr\} \\
& =  \frac{1}{1-e^{-\lambda_1t}} \biggl\{ \frac{\lambda_1 e^{-\lambda_1 t}}{2c_1} \sum_{k=0}^{\infty} \frac{1}{(k!)^2} \; \left(\frac{\lambda_1 t}{2}\right)^{2k} \left( 1+ \frac{\lambda_1 t}{2k+2} \right) \\
& \qquad\qquad\quad \times \left[ c_1t \; F\left( -k, \frac{1}{2}; \; \frac{3}{2}; \; 1 \right) - (c_2t-r) \; F\left( -k, \frac{1}{2}; \; \frac{3}{2}; \; \frac{(c_2t-r)^2}{c_1^2t^2} \right)  \right] \biggr\} , 
\endaligned
\end{equation}
$$\text{if} \; r\in ((c_1-c_2)t, \; (c_1+c_2)t] .$$ 
Formula (\ref{dist8}) can be simplified. In view of (\ref{gauss1}), one can easily show that 
\begin{equation}\label{dddist8}
\frac{\lambda_1 e^{-\lambda_1 t}}{2c_1} \sum_{k=0}^{\infty} \frac{1}{(k!)^2} \; \left(\frac{\lambda_1 t}{2}\right)^{2k} \left( 1+ \frac{\lambda_1 t}{2k+2} \right) \; 
c_1t \; F\left( -k, \frac{1}{2}; \; \frac{3}{2}; \; 1 \right) = \frac{1 - e^{-\lambda_1 t}}{2} 
\end{equation}
and, therefore, (\ref{dist8}) takes the form: 
\begin{equation}\label{ddist8}
\aligned
& \text{Pr} \bigl\{ \rho(t) < r \; \bigl| \; N_1(t)\ge 1, \; N_2(t) = 0 \bigr\} \\
& = \frac{1}{2} - \frac{\lambda_1 (c_2t-r) e^{-\lambda_1 t}}{2c_1(1 - e^{-\lambda_1 t})} \sum_{k=0}^{\infty} \frac{1}{(k!)^2} \left(\frac{\lambda_1 t}{2}\right)^{2k} \left( 1+ \frac{\lambda_1 t}{2k+2} \right) F\left( -k, \frac{1}{2}; \frac{3}{2}; \frac{(c_2t-r)^2}{c_1^2t^2} \right) , 
\endaligned
\end{equation}
$$\text{if} \; r\in ((c_1-c_2)t, \; (c_1+c_2)t] .$$ 

\vskip 0.1cm

$\bullet$ {\it Evaluation of} $\text{Pr} \left\{ \rho(t) < r \; \bigl| \; N_1(t)=0, \; N_2(t)\ge 1 \right\}$. It is obvious that for arbitrary $r\in (0, \; (c_1-c_2)t]$ the 
following relation holds: 
\begin{equation}\label{dist9}
\text{Pr} \bigl\{ \rho(t) < r \; \bigl| \; N_1(t)=0, \; N_2(t) \ge 1 \bigr\} = 0, \qquad \text{if} \; r\in (0, \; (c_1-c_2)t] .
\end{equation}

Let now $r\in ((c_1-c_2)t, \; (c_1+c_2)t]$. Since  
$$\aligned 
\left\{ N_1(t)=0 \right\} & = \left\{ X_1(t)=-c_1t \right\} + \left\{ X_1(t)=c_1t \right\} , \\
\left\{ N_2(t)\ge 1 \right\} & = \left\{ X_2(t) \in (-c_2t, \; c_2t) \right\},
\endaligned$$
then, similarly as above, one can show that 
$$\aligned 
\text{Pr} & \bigl\{ \rho(t) < r \; \bigl| \; N_1(t)=0, \; N_2(t) \ge 1 \bigr\} \\ 
& = \frac{1}{2(1-e^{-\lambda_2t})} \biggl[ \text{Pr} \bigl\{ X_2(t) \in (-c_2t, \; -c_1t+r) \bigr\} + \text{Pr} \bigl\{ X_2(t) \in (c_1t-r, \; c_2t) \bigr\} \biggr] .
\endaligned$$ 
Applying formula (\ref{sub1}) of Proposition 1, we get:  
$$\aligned
\text{Pr} & \bigl\{ \rho(t) < r \; \bigl| \; N_1(t)=0, \; N_2(t) \ge 1 \bigr\} \\
& = \frac{1}{2(1-e^{-\lambda_2t})} \biggl\{ \frac{\lambda_2 e^{-\lambda_2 t}}{2c_2} \sum_{k=0}^{\infty} \frac{1}{(k!)^2} \; \left(\frac{\lambda_2 t}{2}\right)^{2k} \left( 1+ \frac{\lambda_2 t}{2k+2} \right) \\
& \qquad\qquad\qquad \times \left[ (-c_1t+r) \; F\left( -k, \frac{1}{2}; \; \frac{3}{2}; \; \frac{(-c_1t+r)^2}{c_2^2t^2} \right) + c_2t \; F\left( -k, \frac{1}{2}; \; \frac{3}{2}; \; 1 \right) \right] \\
& \qquad\qquad\qquad  + \frac{\lambda_2 e^{-\lambda_2 t}}{2c_2} \sum_{k=0}^{\infty} \frac{1}{(k!)^2} \; \left(\frac{\lambda_2 t}{2}\right)^{2k} \left( 1+ \frac{\lambda_2 t}{2k+2} \right) \\
& \qquad\qquad\qquad \times \left[ c_2t \; F\left( -k, \frac{1}{2}; \; \frac{3}{2}; \; 1 \right) - (c_1t-r) \; F\left( -k, \frac{1}{2}; \; \frac{3}{2}; \; \frac{(c_1t-r)^2}{c_2^2t^2} \right) \right]  \biggr\} \\ 
& = \frac{1}{1-e^{-\lambda_2t}} \biggl\{ \frac{\lambda_2 e^{-\lambda_2 t}}{2c_2} \sum_{k=0}^{\infty} \frac{1}{(k!)^2} \; \left(\frac{\lambda_2 t}{2}\right)^{2k} \left( 1+ \frac{\lambda_2 t}{2k+2} \right) \\
& \qquad\qquad\qquad \times \left[ c_2t \; F\left( -k, \frac{1}{2}; \; \frac{3}{2}; \; 1 \right) - (c_1t-r) \; F\left( -k, \frac{1}{2}; \; \frac{3}{2}; \; \frac{(c_1t-r)^2}{c_2^2t^2} \right) \right] \biggr\} . 
\endaligned$$
Taking into account that 
$$\frac{\lambda_2 e^{-\lambda_2 t}}{2c_2} \sum_{k=0}^{\infty} \frac{1}{(k!)^2} \; \left(\frac{\lambda_2 t}{2}\right)^{2k} \left( 1 + \frac{\lambda_2 t}{2k+2} \right) 
\; c_2t \; F\left( -k, \frac{1}{2}; \; \frac{3}{2}; \; 1 \right) = \frac{1 - e^{-\lambda_2 t}}{2} $$
we finally obtain: 
\begin{equation}\label{dist10}
\aligned 
& \text{Pr} \bigl\{ \rho(t) < r \; \bigl| \; N_1(t)=0, \; N_2(t) \ge 1 \bigr\} \\
& = \frac{1}{2} - \frac{\lambda_2 (c_1t-r) e^{-\lambda_2 t}}{2c_2(1-e^{-\lambda_2t})} \sum_{k=0}^{\infty} \frac{1}{(k!)^2} \left(\frac{\lambda_2 t}{2}\right)^{2k} 
\left( 1 + \frac{\lambda_2 t}{2k+2} \right) F\left( -k, \frac{1}{2}; \frac{3}{2}; \frac{(c_1t-r)^2}{c_2^2t^2} \right) ,
\endaligned
\end{equation}
$$\text{if} \; r\in ((c_1-c_2)t, \; (c_1+c_2)t] .$$

\vskip 0.1cm

$\bullet$ {\it Evaluation of} $\text{Pr} \left\{ \rho(t) < r \; \bigl| \; N_1(t)\ge 1, \; N_2(t)\ge 1 \right\}$. Since 
$$\left\{ N_1(t)\ge 1 \right\} = \left\{ X_1(t)\in (-c_1t, \; c_1t) \right\} , \qquad \left\{ N_2(t)\ge 1 \right\} = \left\{ X_2(t) \in (-c_2t, \; c_2t) \right\},$$
then, for the fourth conditional probability on the right-hand side of (\ref{dist5}), we have: 

$$\aligned 
& \text{Pr} \left\{ \rho(t) < r \; \bigl| \; N_1(t)\ge 1, \; N_2(t)\ge 1 \right\} \\
& = \frac{\text{Pr} \bigl\{ \left\{ \rho(t) < r\right\} \cap \left\{ X_1(t)\in (-c_1t, \; c_1t) \right\} \cap \left\{ X_2(t)\in (-c_2t, \; c_2t) \right\} \bigr\} }{\text{Pr}\left\{ X_1(t)\in (-c_1t, \; c_1t) \right\} \; \text{Pr}\left\{ X_2(t)\in (-c_2t, \; c_2t) \right\}} \\
& = \frac{1}{(1-e^{-\lambda_1t})(1-e^{-\lambda_2t})} \\ 
& \times \text{Pr} \bigl\{ \left\{ X_1(t)\in (X_2(t)-r, \; X_2(t)+r) \right\} \cap \left\{ X_1(t)\in (-c_1t, \; c_1t) \right\} \cap \left\{ X_2(t)\in (-c_2t, \; c_2t) \right\} \bigr\} \\
& = \frac{1}{(1-e^{-\lambda_1t})(1-e^{-\lambda_2t})} \\ 
& \times \text{Pr} \bigl\{ \left\{ X_1(t)\in \left( \max\{ X_2(t)-r, \; -c_1t\}, \; \min\{ X_2(t)+r, \; c_1t\} \right) \right\} \cap \left\{ X_2(t)\in (-c_2t, \; c_2t) \right\} \bigr\} \\
& = \frac{1}{(1-e^{-\lambda_1t})(1-e^{-\lambda_2t})} \int\limits_{-c_2t}^{c_2t} \text{Pr} \bigl\{ X_1(t) \in (\alpha(x,r), \; \beta(x,r)) \; \big| \; X_2(t)=x \bigr\} \; \text{Pr} \bigl\{ X_2(t) \in dx \bigr\} ,  
\endaligned$$ 
where 
$$\alpha(x,r) = \max\{ x-r, \; -c_1t \} , \qquad \beta(x,r) = \min\{ x+r, \; c_1t \} .$$
In view of (\ref{sub1}) and (\ref{prop6}), we obtain: 
\begin{equation}\label{dist11}
\aligned 
\text{Pr} & \bigl\{ \rho(t) < r \; \bigl| \; N_1(t) \ge 1, \; N_2(t) \ge 1 \bigr\} \\
& = \frac{1}{(1-e^{-\lambda_1t})(1-e^{-\lambda_2t})} \int\limits_{-c_2t}^{c_2t} \biggl\{ \frac{\lambda_1 e^{-\lambda_1 t}}{2c_1} \sum_{k=0}^{\infty} \frac{1}{(k!)^2} \; \left( \frac{\lambda_1 t}{2} \right)^{2k} \left( 1 + \frac{\lambda_1 t}{2k+2} \right) \\
& \quad \times \left[ \beta(x,r) F\left( -k, \frac{1}{2}; \frac{3}{2}; \frac{(\beta(x,r))^2}{c_1^2t^2} \right) - \alpha(x,r) F\left( -k, \frac{1}{2}; \frac{3}{2}; \frac{(\alpha(x,r))^2}{c_1^2t^2} \right) \right] \biggr\} f_2^{ac}(x,t) dx \\ 
& = \frac{e^{-\lambda_1t}}{(1-e^{-\lambda_1t})(1-e^{-\lambda_2t})} \frac{\lambda_1}{2c_1} \sum_{k=0}^{\infty} \frac{1}{(k!)^2} \; \left( \frac{\lambda_1 t}{2} \right)^{2k} 
\left( 1 + \frac{\lambda_1 t}{2k+2} \right) \mathcal I_k(r,t) ,
\endaligned
\end{equation}
where the integral factor $\mathcal I_k(r,t)$ is defined by (\ref{distINT}) and $f_2^{ac}(x,t)$ is the density of the absolutely continuous part of the distribution of the telegraph process $X_2(t)$ given by (\ref{prop6}).

Substituting now (\ref{ddist5}), (\ref{dist7}), (\ref{dist9}) and (\ref{dist11}) into (\ref{dist5}) we obtain the term $G(r,t)$ in distribution function (\ref{dist4}) defined in the interval $r\in (0, \; (c_1-c_2)t]$ and given by formula (\ref{distG}). 
Similarly, by substituting (\ref{ddist5}), (\ref{ddist8}), (\ref{dist10}) and (\ref{dist11}) into (\ref{dist5}) we obtain the term $Q(r,t)$ in distribution function (\ref{dist4}) defined in the interval $r\in ((c_1-c_2)t, \; (c_1+c_2)t]$ and given by formula (\ref{distQ}).
The theorem is thus completely proved. $\square$

\bigskip 

{\it Remark 4.} One can easily see that if $r\in (0, \; (c_1-c_2)t]$, then the variables $\alpha(x,r)$ and $\beta(x,r)$ take the values 
$$\alpha(x,r) = x-r, \qquad \beta(x,r) = x+r, \qquad \text{for} \; r\in (0, \; (c_1-c_2)t],$$
independently of $x\in (-c_2t, \; c_2t)$. In this case the integral factor $\mathcal I_k(r,t)$ can, therefore, be rewritten in a bit more explicit form. In contrast, if $r\in ((c_1-c_2)t, \; (c_1+c_2)t]$, then each of these variables can take both possible values.  

\vskip 0.2cm 

{\it Remark 5.} Taking into account that, for any $x\in (-c_2t, \; c_2t)$,  
$$\alpha(x,0)=\beta(x,0)=x, \qquad \alpha(x,(c_1+c_2)t) = -c_1t, \qquad \beta(x,(c_1+c_2)t) = c_1t,$$
$$\alpha(x,(c_1-c_2)t) = x - (c_1-c_2)t, \qquad \beta(x,(c_1-c_2)t) = x + (c_1-c_2)t,$$
one can easily prove the following limiting relations: 
\begin{equation}\label{distLIM}
\aligned 
& \lim\limits_{r\to 0+0} G(r,t) = 0, \qquad \lim\limits_{r\to (c_1+c_2)t-0} Q(r,t) = 1 - \frac{1}{2} e^{-(\lambda_1+\lambda_2)t}, \\
& \lim\limits_{r\to (c_1-c_2)t+0} Q(r,t) - \lim\limits_{r\to (c_1-c_2)t-0} G(r,t) = \frac{1}{2} e^{-(\lambda_1+\lambda_2)t} . 
\endaligned
\end{equation}
Formulas (\ref{distLIM}) show that distribution function (\ref{dist4}) is left-continuous with jumps of the same amplitude $e^{-(\lambda_1+\lambda_2)t}/2$ 
at the singularity points $(c_1\pm c_2)t$. This entirely accords with the structure of the distribution of the process $\rho(t)$ described above. 

\vskip 0.2cm 

{\it Remark 6.} In the use of distribution function (\ref{dist4}) the crucial point is the possibility of computing the integral term $\mathcal I_k(r,t)$ given by (\ref{distINT}). By means of tedious computations and by applying formulas (\ref{app6}) and (\ref{app7}) of the Appendix one can obtain a series representations of integral $\mathcal I_k(r,t)$, however it has an extremely complicated and cumbersome form and therefore is omitted here. That is why for practical purposes it is more convenient to use just the integral form of factor $\mathcal I_k(r,t)$, which is easily computable by usual personal computer (for more details see Section 5 below dealing with some numerical analysis of formulas obtained).  

\vskip 0.2cm 

We conclude this section by presenting a result related to the more simple case of equal velocities. Suppose that both the telegraph processes $X_1(t)$ and $X_2(t)$ develop with the same speed $c_1=c_2=c$. In this case the support of distribution (\ref{distPhi}) is the close interval $[0, \; 2ct]$. The singular component of distribution has the density (as generalized function) 
\begin{equation}\label{dist12}
\varphi^s(r,t) = \frac{e^{-(\lambda_1+\lambda_2)t}}{2} \; \bigl[ \delta(r) + \delta(r-2ct) \bigr] , \qquad r\in\Bbb R, \quad t>0, 
\end{equation}
concentrated at the terminal points 0 and $2ct$, while the open interval $(0, \; 2ct)$ is the support of the absolutely continuous part of distribution (\ref{distPhi}). 
The form of distribution function (\ref{distPhi}) for the case of equal velocities is presented by the following theorem.

\bigskip 

{\bf Theorem 2.} {\it Under the condition $c_1=c_2=c$ the probability distribution function} (\ref{distPhi}) {\it has the form:}
\begin{equation}\label{dist13}
\Phi(r,t) = \left\{ \aligned 0 , \qquad & \text{if} \; r\in (-\infty, \; 0] ,\\ 
                             H(r,t), \qquad & \text{if} \; r\in (0, \; 2ct] ,\\
                             1 , \qquad & \text{if} \; r\in (2ct, \; +\infty) ,
\endaligned \right. \qquad r\in\Bbb R, \qquad t>0,
\end{equation}
{\it where function $H(r,t)$ is given by the formula:}

\begin{equation}\label{distH}
\aligned 
& H(r,t) \\
& = \frac{1}{2} \left[ \left( 1-e^{-\lambda_1t} \right) e^{-\lambda_2t} + e^{-\lambda_1t} \left( 1-e^{-\lambda_2t} \right) + e^{-(\lambda_1+\lambda_2)t} \right] \\ 
& \; - e^{-(\lambda_1+\lambda_2)t} \; \left(  1 - \frac{r}{ct} \right) \sum_{k=0}^{\infty} \frac{1}{(k!)^2} \left[ \left( \frac{\lambda_1 t}{2} \right)^{2k+1} \left( 1 + \frac{\lambda_1 t}{2k+2} \right) + \left( \frac{\lambda_2 t}{2} \right)^{2k+1} \left( 1 + \frac{\lambda_2 t}{2k+2} \right) \right] \\ 
& \hskip 7cm \times F\left( -k, \frac{1}{2}; \; \frac{3}{2}; \; \left(  1 - \frac{r}{ct} \right)^2 \right) \\ 
& \; + e^{-\lambda_1t} \; \frac{\lambda_1}{2c} \sum_{k=0}^{\infty} \frac{1}{(k!)^2} \; \left( \frac{\lambda_1 t}{2} \right)^{2k} \left( 1 + \frac{\lambda_1 t}{2k+2} \right) \mathcal J_k(r,t) , 
\endaligned
\end{equation}
{\it with the integral term} 
\begin{equation}\label{distINT1}
\aligned 
\mathcal J_k(r,t) = & \frac{\lambda_2 e^{-\lambda_2 t}}{2c} \int\limits_{-ct}^{ct} \biggl[ \beta(x,r) F\left( -k, \frac{1}{2}; \frac{3}{2}; \frac{(\beta(x,r))^2}{c^2t^2} \right) - \alpha(x,r) F\left( -k, \frac{1}{2}; \frac{3}{2}; \frac{(\alpha(x,r))^2}{c^2t^2} \right) \biggr] \\ 
& \hskip 3cm \times \biggl[ I_0\left( \frac{\lambda_2}{c} \sqrt{c^2t^2-x^2} \right) + \frac{ct}{\sqrt{c^2t^2-x^2}} I_1\left( \frac{\lambda_2}{c} 
\sqrt{c^2t^2-x^2} \right) \biggr] dx ,
\endaligned
\end{equation}
{\it where the variables $\alpha(x,r)$ and $\beta(x,r)$ are defined by the formulas}: 
\begin{equation}\label{distAB1}
\alpha(x,r) = \max\{ -ct, \; x-r \} , \qquad \beta(x,r) = \min\{ ct, \; x+r \} , 
\end{equation}
$$x\in(-ct, \; ct), \qquad r\in (0, \; 2ct).$$

\vskip 0.2cm

{\it Proof.} The proof is similar to that of Theorem 1 and therefore is omitted here. $\square$

\section{Some Numerical Results}

\numberwithin{equation}{section}

While distribution function (\ref{dist4}) has fairly complicated analytical form, it can, however, be approximately evaluated with good accuracy by using standard mathematical programs (such as MATHEMATICA or MAPLE) and usual personal computer. As is noted in Remark 6 above, the crucial point is the evaluation of the integral term $\mathcal I_k(r,t)$ represented  by formula (\ref{distINT}) (for $c_1>c_2$) or by formula (\ref{distINT1}) (for $c_1=c_2=c$). 

To approximately evaluate the series in functions (\ref{distG}) and (\ref{distQ}), we do not need to compute integral term (\ref{distINT}) for all $k\ge 0$. We notice that each series contains the factor $1/(k!)^2$ providing its very fast convergence. Taking into account that all the hypergeometric functions in (\ref{distG}) and (\ref{distQ}) are uniformly bounded for all $k$ (see estimate (\ref{aapp1}) below), we can conclude that each series in (\ref{distG}) and (\ref{distQ}) converges even faster than that of the modified Bessel function $I_0(z)$ given by (\ref{pprop5}). In fact, one can see that, if we take only five terms of each series in functions $G(r,t)$ and $Q(r,t)$, its approximate value becomes stabilized at fourth digit already.   

The results of numerical analysis for distribution function (\ref{dist4}) with parameters: 
$$\lambda_1=2, \quad \lambda_2=1, \quad c_1=4, \quad c_2=2, \quad t=3,$$
with function $G(r,3)$ defined in the subinterval $r\in (0, \; 6]$ and function $Q(r,3)$ defined in the subinterval $r\in (6, \; 12]$, respectively, are given in Tables 1 and 2 below:

\begin{center} 
{\bf Table 1:} Function $G(r,3)$ in the subinterval $r\in (0, \; 6]$ 
\end{center}

\begin{center}

\begin{tabular}{|r|r||r|r||r|r|}
\hline
 $r$ & $G(r,3)$ & $r$ & $G(r,3)$ & $r$ & $G(r,3)$ \\
\hline\hline
  0.2 & 0.0271 & 2.2 & 0.2916 & 4.2 & 0.5269 \\
\hline
  0.4 & 0.0541 & 2.4 & 0.3168 & 4.4 & 0.5480 \\
\hline
  0.6 & 0.0811 & 2.6 & 0.3417 & 4.6 & 0.5686 \\
\hline
  0.8 & 0.1080 & 2.8 & 0.3663 & 4.8 & 0.5888 \\
\hline
  1.0 & 0.1348 & 3.0 & 0.3905 & 5.0 & 0.6083 \\
\hline
  1.2 & 0.1614 & 3.2 & 0.4143 & 5.2 & 0.6274 \\
\hline
  1.4 & 0.1879 & 3.4 & 0.4377 & 5.4 & 0.6459 \\
\hline
  1.6 & 0.2142 & 3.6 & 0.4607 & 5.6 & 0.6639 \\
\hline
  1.8 & 0.2402 & 3.8 & 0.4832 & 5.8 & 0.6813 \\
\hline
  2.0 & 0.2660 & 4.0 & 0.5053 & 6.0 & 0.6982 \\
\hline
\end{tabular}

\end{center}

\vskip 1cm

\begin{center} 
{\bf Table 2:} Function $Q(r,3)$ in the subinterval $r\in (6, \; 12]$ 
\end{center}

\begin{center}

\begin{tabular}{|r|r||r|r||r|r|}
\hline
 $r$ & $Q(r,3)$ & $r$ & $Q(r,3)$ & $r$ & $Q(r,3)$ \\
\hline\hline
  6.2 & 0.7146 & 8.2 & 0.8472 & 10.2 & 0.9294 \\
\hline
  6.4 & 0.7302 & 8.4 & 0.8576 & 10.4 & 0.9350 \\
\hline
  6.6 & 0.7455 & 8.6 & 0.8674 & 10.6 & 0.9405 \\
\hline
  6.8 & 0.7601 & 8.8 & 0.8768 & 10.8 & 0.9461 \\
\hline
  7.0 & 0.7741 & 9.0 & 0.8855 & 11.0 & 0.9510 \\
\hline
  7.2 & 0.7877 & 9.2 & 0.8945 & 11.2 & 0.9554 \\
\hline
  7.4 & 0.8006 & 9.4 & 0.9019 & 11.4 & 0.9589 \\
\hline
  7.6 & 0.8131 & 9.6 & 0.9093 & 11.6 & 0.9631 \\
\hline
  7.8 & 0.8250 & 9.8 & 0.9170 & 11.8 & 0.9673 \\
\hline
  8.0 & 0.8364 & 10.0 & 0.9233 & 12.0 & 0.9704 \\
\hline
\end{tabular}

\end{center}

Note that in evaluating these functions we take only seven terms of all the series. Also note that, while function $Q(r,3)$ is defined  
in the whole interval $(6, \; 18]$, we consider it only in the first half $(6, \; 12]$ because in the second half $(12, \; 18]$ function 
$Q(r,3)$ has negligible increments. Notice that, at the terminal point $r=18$, it takes the value $Q(18,3) \approx 0.999937$ 
(with ten terms of each series taken) which differs from 1 on the value of jump amplitude at this point.

\section{Appendix}

\numberwithin{equation}{section}

In this appendix we prove two auxiliary lemmas concerning some indefinite integrals and an useful formula related to conditional probabilities. 

\bigskip

{\bf Lemma A1.} {\it For arbitrary $q \ge 0, \; p>0$ such that $|x| \le p$ the following formulas hold:}
\begin{equation}\label{app1}
\int x^n \; I_0(q\sqrt{p^2-x^2}) \; dx  
= \frac{x^{n+1}}{n+1} \sum_{k=0}^{\infty} \frac{1}{(k!)^2} \left( \frac{pq}{2} \right)^{2k} F\left( -k, \frac{n+1}{2}; \frac{n+3}{2}; \frac{x^2}{p^2} \right) + \psi_1, 
\end{equation}
$$n\ge 0, \qquad |x| \le p,$$
\begin{equation}\label{app2}
\int x^n \; \frac{I_1(q\sqrt{p^2-x^2})}{\sqrt{p^2-x^2}} \; dx = \frac{x^{n+1}}{p(n+1)} \sum_{k=0}^{\infty} \frac{1}{k! \; (k+1)!} \left( \frac{pq}{2} \right)^{2k+1} F\left( -k, \frac{n+1}{2}; \frac{n+3}{2}; \frac{x^2}{p^2} \right) + \psi_2, 
\end{equation}
$$n\ge 0, \qquad |x| \le p,$$
{\it where $\psi_1, \; \psi_2$ are arbitrary functions not depending on $x$.} 

\vskip 0.2cm

{\it Proof.} Let us check formula (\ref{app1}). First, we prove that the series on the right-hand side of (\ref{app1}) converges uniformly with respect 
to $x\in [-p, p]$. In view of the following uniform (in $z\in [0, 1]$) estimate 
\begin{equation}\label{aapp1}
0 < F\left( -k, \frac{n+1}{2}; \frac{n+3}{2}; z  \right) \le 1 , \qquad 0\le z\le 1, \quad n\ge 0, \quad k\ge 0,
\end{equation} 
we obtain the inequality  
$$\sum_{k=0}^{\infty} \frac{1}{(k!)^2} \left( \frac{pq}{2} \right)^{2k} F\left( -k, \frac{n+1}{2}; \frac{n+3}{2}; \frac{x^2}{p^2} \right) 
\le \sum_{k=0}^{\infty} \frac{1}{(k!)^2} \left( \frac{pq}{2} \right)^{2k} = I_0(pq) < \infty , $$
proving the uniform convergence in $x\in [-p, p]$ of the series. From this fact it follows that one may differentiate each term of the series on the right-hand side of (\ref{app1}) separately. Thus, differentiating in $x$ the expression on the right-hand side of (\ref{app1}) and taking into account that 
$$\aligned 
(-k)_s & = \frac{(-1)^s \; k!}{(k-s)!} , \qquad k\ge 0, \quad 0\le s\le k, \\ 
\frac{(a)_s}{(a+1)_s} & = \frac{a}{a+s} , \qquad a>0, \quad s\ge 0,  
\endaligned$$
we obtain 
$$\aligned 
\frac{d}{dx} & \left[ \frac{x^{n+1}}{n+1} \sum_{k=0}^{\infty} \frac{1}{(k!)^2} \left( \frac{pq}{2} \right)^{2k} F\left( -k, \frac{n+1}{2}; \frac{n+3}{2}; \frac{x^2}{p^2} \right) \right] \\ 
& = \frac{1}{n+1} \sum_{k=0}^{\infty} \frac{1}{(k!)^2} \left( \frac{pq}{2} \right)^{2k} \frac{d}{dx} \left[ x^{n+1} F\left( -k, \frac{n+1}{2}; \frac{n+3}{2}; \frac{x^2}{p^2} \right) \right] \\ 
& = \frac{1}{n+1} \sum_{k=0}^{\infty} \frac{1}{(k!)^2} \left( \frac{pq}{2} \right)^{2k} \; \frac{d}{dx} \left[ \sum_{s=0}^k \frac{(-1)^s \; k! \; (n+1)}{s! \; (k-s)! \; (n+2s+1)} \; \frac{x^{n+2s+1}}{p^{2s}} \right] \\ 
& = x^n \sum_{k=0}^{\infty} \frac{1}{(k!)^2} \left( \frac{pq}{2} \right)^{2k} \; \left[ \sum_{s=0}^k (-1)^s \binom ks \left( \frac{x^2}{p^2} \right)^{s} \right] \\
& = x^n \sum_{k=0}^{\infty} \frac{1}{(k!)^2} \left( \frac{pq}{2} \right)^{2k} \; \left( 1 - \frac{x^2}{p^2} \right)^k \\
\endaligned$$
$$\aligned 
& = x^n \sum_{k=0}^{\infty} \frac{1}{(k!)^2} \left( \frac{q}{2} \right)^{2k} (p^2 - x^2)^k \\
& = x^n \sum_{k=0}^{\infty} \frac{1}{(k!)^2} \left( \frac{q \sqrt{p^2 - x^2}}{2} \right)^{2k} \\
& = x^n \; I_0(q\sqrt{p^2-x^2}) , 
\endaligned$$
yielding the integrand on the left-hand side of (\ref{app1}). Formula (\ref{app2}) can be checked similarly. The lemma is proved. $\square$ 

\bigskip 

In particular, by setting $n=0$ in (\ref{app1}) and (\ref{app2}), we arrive to the formulas:  
\begin{equation}\label{app4}
\int I_0(q\sqrt{p^2-x^2}) \; dx  
= x \sum_{k=0}^{\infty} \frac{1}{(k!)^2} \left( \frac{pq}{2} \right)^{2k} F\left( -k, \frac{1}{2}; \frac{3}{2}; \frac{x^2}{p^2} \right) + \psi_1, \qquad |x| \le p, 
\end{equation}
\begin{equation}\label{app5}
\int \frac{I_1(q\sqrt{p^2-x^2})}{\sqrt{p^2-x^2}} \; dx = \frac{x}{p} \sum_{k=0}^{\infty} \frac{1}{k! \; (k+1)!} \left( \frac{pq}{2} \right)^{2k+1} F\left( -k, \frac{1}{2}; \frac{3}{2}; \frac{x^2}{p^2} \right) + \psi_2, \qquad |x| \le p.
\end{equation}

Applying Lemma A1 we have for arbitrary real $a$:
\begin{equation}\label{app6}
\aligned 
\int & (a\pm x)^n \; I_0(q\sqrt{p^2-x^2}) \; dx \\
& = \sum_{m=0}^n (\pm 1)^m \; \binom nm a^{n-m} \; \frac{x^{m+1}}{m+1} \sum_{k=0}^{\infty} \frac{1}{(k!)^2} \left( \frac{pq}{2} \right)^{2k} F\left( -k, \frac{m+1}{2}; \frac{m+3}{2}; \frac{x^2}{p^2} \right) + \psi_1 , 
\endaligned
\end{equation}
$$n\ge 0, \qquad |x| \le p.$$ 

Similarly, for arbitrary real $a$ the following formula holds: 
\begin{equation}\label{app7}
\aligned 
& \int (a\pm x)^n \; \frac{I_1(q\sqrt{p^2-x^2})}{\sqrt{p^2-x^2}} \; dx \\
& = \frac{1}{p} \sum_{m=0}^n (\pm 1)^m \binom{n}{m} a^{n-m} \; 
\frac{x^{m+1}}{m+1} \sum_{k=0}^{\infty} \frac{1}{k! \; (k+1)!} \left( \frac{pq}{2} \right)^{2k+1} F\left( -k, \frac{m+1}{2}; \frac{m+3}{2}; \frac{x^2}{p^2} \right) + \psi_2, 
\endaligned 
\end{equation}
$$n\ge 0, \qquad |x| \le p.$$ 

\bigskip

In the next lemma we prove an useful formula for the probabilities conditioned by pairwisely independent random events that has been used in our analysis.  

\bigskip

{\bf Lemma A2.} {\it Let $(\Omega,\mathcal F,\bold P)$ be a probability space and let $A,B,C,D\in\mathcal F$ be the random events such that $B,C,D$ are pairwisely independent,  $C\cap D=\varnothing, \; \bold P(C)=\bold P(D)\neq 0, \; \bold P(B)\neq 0$. Then the following formula holds:}
\begin{equation}\label{app8}
\bold P(A\; | \; B(C+D)) = \frac{1}{2} \bigl[ \bold P(A\; | \; BC) + \bold P(A\; | \; BD) \bigr] . 
\end{equation}

\vskip 0.2cm

{\it Proof.} Under the lemma's conditions, we have:  

$$\aligned
\bold P(A\; | \; B(C+D)) & = \frac{\bold P(AB(C+D))}{\bold P(B(C+D))} \\
& = \frac{\bold P(ABC) + \bold P(ABD)}{\bold P(BC)+\bold P(BD)} \\
& = \frac{\bold P(A \; | \; BC) \bold P(B) \bold P(C) + \bold P(A \; | \; BD) \bold P(B) \bold P(D)}{\bold P(B) \left[ \bold P(C)+\bold P(D) \right]} \\
& = \frac{\bold P(A \; | \; BC) \bold P(B) \bold P(C) + \bold P(A \; | \; BD) \bold P(B) \bold P(C)}{2 \; \bold P(B) \bold P(C)} \\
& = \frac{1}{2} \bigl[ \bold P(A\; | \; BC) + \bold P(A\; | \; BD) \bigr] .
\endaligned$$
The lemma is proved. $\square$

\bigskip

{\bf Acknowledgements.} This article was written in the framework of the bilateral Germany-Moldova 
research project 13.820.18.01/GA. 

\bigskip

\numberwithin{equation}{section}

\end{document}